\documentclass[letterpaper, 10pt]{article}
\usepackage{amsmath, amsthm, amssymb, amsfonts}

\newtheorem{thm}{Theorem}[section]
\newtheorem{cor}[thm]{Corollary}

\newtheorem{lem}[thm]{Lemma}
\newtheorem{defn}[thm]{Definition}
\newtheorem{prop}[thm]{Proposition}

\newcommand{\SIMPLE}{\textsc{Simple}}
\newcommand{\cExp}[3]{\mathbf{E}_{#1}\left[\left.#2\;\right\vert #3\right]}
\newcommand{\sExp}[2]{\mathbf{E}_{#1}\left[#2\right]}
\newcommand{\Exp}[1]{\sExp{}{#1}}
\newcommand{\Prob}[1]{\Pr\left[#1\right]}
\newcommand{\cProb}[3]{\Pr_{#1}\left[ \left. #2 \;\right\vert #3 \right]}

\newcommand{\Gnd}{\mathcal{G}_{n,d}}
\newcommand{\Gnp}{\mathcal{G}(n,p)}
\newcommand{\Pnd}{\mathcal{P}_{n,d}}
\newcommand{\mC}{\mathcal{C}}
\newcommand{\mD}{\mathcal{E}}
\newcommand{\mE}{\mathcal{D}}
\newcommand{\mF}{\mathcal{F}}

\title{Random regular graphs of non-constant degree:\\Concentration of the chromatic number}
\author{\sc Sonny Ben-Shimon
        \thanks{School of Computer Science, Raymond and Beverly Sackler Faculty of Exact Sciences, Tel Aviv University, Tel Aviv, 69978,
        Israel. Email: {\tt sonny@tau.ac.il}.}
        \and
        \sc Michael Krivelevich
        \thanks{School of Mathematical Sciences, Raymond and Beverly Sackler Faculty of Exact Sciences, Tel Aviv University, Tel Aviv, 69978,
        Israel. Email: {\tt krivelev@tau.ac.il}.
        Research supported in part by a USA-Israel BSF grant, by a grant from the Israel Science Foundation, and by Pazy Memorial Award.}
}

\begin{document}
\maketitle
\begin{abstract}
In this work we show that with high probability the chromatic number
of a graph sampled from the random regular graph model $\Gnd$ for
$d=o(n^{1/5})$ is concentrated in two consecutive values, thus
extending a previous result of Achlioptas and Moore. This
concentration phenomena is very similar to that of the binomial
random graph model $\Gnp$ with $p=\frac{d}{n}$. Our proof is largely
based on ideas of Alon and Krivelevich who proved this two-point
concentration result for $\Gnp$ for $p=n^{-\delta}$ where
$\delta>1/2$. The main tool used to derive such a result is a
careful analysis of the distribution of edges in $\Gnd$, relying
both on the switching technique and on bounding the probability of
exponentially small events in the configuration model.
\end{abstract}

\section{Introduction}
The most widely used random graph model is the binomial random
graph, $\Gnp$. In this model, which was introduced in a slightly
modified form by Erd\H{o}s and R\`{e}nyi, we start with $n$
vertices, labeled, say, by $\{1,\ldots,n\}=[n]$, and select a graph
on these $n$ vertices by going over all ${n\choose 2}$ pairs of
vertices, deciding uniformly at random with probability $p$ for a
pair to be an edge. $\Gnp$ is thus a probability space of all
labeled graphs on the vertex set $[n]$ where the probability of such
a graph, $G=([n],E)$, to be selected is $p^{|E|}(1-p)^{{n\choose
2}-|E|}$. This product probability space gives us a wide variety of
probabilistic tools to analyze the behavior of various random graph
properties of this probability space. (See monographs \cite{Bol2001}
and \cite{JanLucRuc2000} for a thorough introduction to the subject
of random graphs).

In this paper, we consider a different random graph model. Our
probability space, which is denoted by $\Gnd$ (where $dn$ is even),
is the uniform space of all $d$-regular graphs on $n$ vertices
labeled by the set $[n]$. In this model, one cannot apply the
techniques used to study $\Gnp$ as these two models do not share the
same probabilistic properties. Whereas the appearances of edges in
$\Gnp$ are independent, the appearances of edges in $\Gnd$ are not.
Nevertheless, many results obtained thus far for the random regular
graph model $\Gnd$ are in some sense equivalent to the results
obtained in $\Gnp$ with suitable expected degrees, namely, $d=np$
(see, e.g. \cite{Wormald99} and \cite{KrivelevichEtAl2001} for a
collection of results). This relation between the two random graph
models was partially formalized in \cite{KimVu2004}. The interested
reader is referred to \cite{Wormald99} for a thorough survey of the
random regular graph model $\Gnd$.

The main research interest in random graph models is the asymptotic
behavior of properties as we let the number of vertices of our graph
grow to infinity. We say that an event $\mathcal{A}$ in our
probability space occurs \emph{with high probability} (or w.h.p. for
brevity) if $\Prob{\mathcal{A}}\rightarrow 1$ as $n$ goes to
infinity. Therefore, from now on and throughout the rest of this
work, we will always assume $n$ to be large enough. We use the usual
asymptotic notation, that is, for two functions of $n$, $f(n)$ and
$g(n)$, we write $f=O(g)$ if there exists a constant $C>0$ such that
$f(n)\leq C\cdot g(n)$ for large enough values of $n$; $f=o(g)$ if
$f/g\rightarrow 0$ as $n$ goes to infinity; $f=\Omega(g)$ if
$g=O(f)$; $f=\omega(g)$ if $g=o(f)$; $f=\Theta(g)$ if both $f=O(g)$
and $f=\Omega(g)$.

As far as notation goes, we will always assume, unless specified
otherwise, that the set of vertices of our graph is $[n]$. We use
the usual notation of $N_G(U)$ for the set of neighbors of a vertex
set $U$ in a graph $G$, that is, $N_G(U)=\{v\in V(G)\setminus
U\;:\;\exists u\in U\; \{v,u\}\in E(G)\}$. For a single vertex $v$,
we abuse slightly this notation by writing $N_G(v)$ for
$N_G(\{v\})$. We denote the degree of a vertex $v$ in a graph $G$ by
$d_G(v)$, namely, $d_G(v)=|N_G(v)|$. The set of edges spanned by a
set of vertices $U$, or between two disjoint sets, $U$ and $W$, is
denoted by $E(U)$ and $E(U,W)$, respectively, and the cardinalities
of these sets are denoted by $e(U)$ and $e(U,W)$, respectively. We
use the notation $\Gnp$ or $\Gnd$ to denote both the corresponding
probability space, and a random graph generated in this probability
space, where the actual meaning is clear from the context.

A \emph{vertex coloring} of a graph $G$ is an assignment of a color
to each of its vertices. The coloring is \emph{proper} if no two
adjacent vertices are assigned the same color. The \emph{chromatic
number} of the graph $G$, denoted by $\chi(G)$ is the minimum number
of colors used in a proper coloring of it. The chromatic number of a
graph is one of the most widely researched graph parameters. A major
result of Bollob\'{a}s \cite{Bol88} that was later extended by
{\L}uczak \cite {Luc91b} showed that w.h.p.
$\chi(\Gnp)=(1+o(1))\frac{n\ln\frac{1}{1-p}}{2\ln(np)}$, where
$p\leq c$ and $np\rightarrow \infty$. Frieze and {\L}uczak in
\cite{FriezeLuczak92} proved a similar result for $\Gnd$.
\begin{thm}[Frieze and {\L}uczak \cite{FriezeLuczak92}]\label{t:FL92}
For any $0<\delta<\frac{1}{3}$ w.h.p.
\begin{equation}\label{e:FriezeLuczak92}
\left|\chi(\Gnd)-\frac{d}{2\ln d}-\frac{8d\ln\ln
d}{\ln^2d}\right|\leq\frac{8d\ln\ln d}{\ln^2d},
\end{equation}
where $d_0<d<n^{\delta}$ for some fixed positive constant $d_0$.
\end{thm}

Krivelevich et al. \cite{KrivelevichEtAl2001} and Cooper et al.
\cite{CooperEtAl2002a} extended the range of $d$ for which
(\ref{e:FriezeLuczak92}) holds.

We say that a random variable $X$ in some discrete probability space
$\Omega$ is \emph{highly concentrated} if for every $\varepsilon>0$
it takes one of a finite set of values (not depending on the
cardinality of $\Omega$ which we think of as tending to infinity)
with probability at least $1-\varepsilon$. It has been shown that in
the binomial random graph model $\Gnp$ the chromatic number is
highly concentrated when the graphs are sparse enough. A series of
papers, starting with the seminal work of Shamir and Spencer
\cite{ShaSpe87}, and succeeded by {\L}uczak \cite{Luc91}, and Alon
and Krivelevich \cite{AloKri97} prove that $\chi(\Gnp)$ is
concentrated w.h.p. in two consecutive values for
$p=n^{-1/2-\varepsilon}$ where $\varepsilon$ is any positive real.
Achlioptas and Naor \cite{AchNao2005} went even further to compute
the two values on which $\chi(\Gnp)$ is concentrated for
$p=\frac{c}{n}$ for a constant $c$. Finally, Coja-Oghlan, Panagiotou
and Steger \cite{CojPanSte2008}, building upon the foundations of
\cite{AchNao2005}, computed three consecutive values on which
$\chi(\Gnp)$ is concentrated for $p=n^{-3/4-\varepsilon}$ where
$\varepsilon$ is a positive constant.

We prove a similar concentration result in the random regular graph
model $\Gnd$. In the course of the proof in \cite{AloKri97}, the
authors prove that in $\Gnp$ subsets of vertices that are not ``too
large" cannot be ``too dense" with high probability (where ``too
large" and ``too dense" are quantified as functions of $n$ and $p$).
We follow this recipe combined with some structural results on the
``typical'' random regular graphs, and specifically, the number of
edges spanned by a single vertex set of an appropriate size in
$\Gnd$, to prove the following concentration result.
\begin{thm}\label{t:main2}
For every positive constant $\varepsilon$ there exists an integer
$n_0=n_0(\varepsilon)$ such that for every $n>n_0$ and $d=o(n^{1/5})$
there exists an integer $t=t(n,d,\varepsilon)$ such that
\begin{equation*}
\Prob{\chi(\Gnd)\in\{t,t+1\}}\geq 1-\varepsilon.
\end{equation*}
\end{thm}
In other words, Theorem \ref{t:main2} states that for large enough
values of $n$, the chromatic number of $\Gnd$ for every
$d=o(n^{1/5})$ w.h.p. takes one of two consecutive values. This
result extends a previous result by Achlioptas and Moore
\cite{AchMoo2004}, who prove the same concentration result for a
smaller range of values of $d=d(n)$.

The rest of this paper is organized as follows. In Section
\ref{s:preliminaries} we give a brief introduction to random regular
graphs and state some known results which will be of use in the
succeeding sections. In Section \ref{s:structureGnd} we perform a
somewhat technical analysis of some structural properties of $\Gnd$,
and in particular analyze the distribution of edges for various
orders of subsets of vertices and ranges of $d$. We then utilize
these results to give a proof of Theorem \ref{t:main2} in Section
\ref{s:ChromConcentrarion}.

\section{Preliminaries}\label{s:preliminaries} We start by analyzing
and exploring the setting of random regular graphs and the
techniques that we have to tackle problems in this probability
space.
\subsection{The Configuration Model}\label{s:ConfModel}
One of the major obstacles posed by the random regular graph model
is the lack of a random generation process of all $d$-regular graphs
for some given value of $d$. The following generation process,
called the \emph{Configuration Model} was introduced by Bollob\'{a}s
in \cite{Bol80}, and implicitly by Bender and Canfield in
\cite{BenCan78}. Consider a set of $dn$ elements (assuming $dn$ is
even), $\{e_1,...,e_{dn}\}$. Let this set be partitioned into $n$
\emph{cells}, $c_i=\{e_j\;|\;d(i-1)+1\leq j \leq di\}$, where $1\leq
i\leq n$. A perfect matching of the elements into $\frac{dn}{2}$
pairs is called a \emph{pairing}. We denote by $\Pnd$ the uniform
probability space of these
$(dn)!!=\frac{(dn)!}{\left(\frac{dn}{2}\right)!2^{dn/2}}$ possible
pairings. Let $P\in \Pnd$, and $e_a$ be some element. We denote by
$e_a^P$ the element that is paired with $e_a$ in the pairing $P$,
that is if $\{e_a,e_b\}\in P$ then $e_a^P=e_b$ and $e_b^P=e_a$. We
define a multigraph $G(P)$, where $V(G(P))=[n]$, and for every pair
$\{e_a,e_b\}\in P$ where $e_a\in c_i$ and $e_b\in c_j$ we add and
edge connecting $i$ to $j$ in $G(P)$. For a general pairing $P$,
$G(P)$ can obviously have loops and multiple edges. We can define a
random $d$-regular multigraph model that assigns to each such
multigraph $G$ the accumulated probability of all pairings $P$ from
$\Pnd$ such that $G(P)=G$. Although this probability space is not of
simple $d$-regular graphs, but of $d$-regular multigraphs, it can be
easily shown (see e.g. \cite{Wormald99}) that all $d$-regular simple
graphs on $n$ vertices are equiprobable in this space. Now, one can
generate a regular graph in $\Gnd$ by sequentially generating random
pairings $P\in\Pnd$ and taking the first $G(P)$ that is a simple
graph.

We define the event {{\SIMPLE}} as the event that the pairing
generated in $\Pnd$ corresponds to a simple graph. By the uniformity
of the two models, it follows that for any event $\mathcal{A}$ in
$\Gnd$, and any event $\mathcal{B}$ in $\Pnd$ where $P\in
\mathcal{B}\cap{\SIMPLE}\Leftrightarrow G(P)\in\mathcal{A}$ we have
\begin{equation}\label{e:ConfSimple}
\Prob{\mathcal{A}}=\cProb{}{\mathcal{B}}{{\SIMPLE}}\leq\frac{\Prob{\mathcal{B}}}{\Prob{{\SIMPLE}}}.
\end{equation}
McKay and Wormald \cite{McKayWormald91} managed to compute the
probability of {\SIMPLE} for $d=o(\sqrt n)$.
\begin{thm}[McKay and Wormald \cite{McKayWormald91}]\label{t:BoundProbSimple} For
$d=o(\sqrt{n})$,
\begin{equation}
\Prob{{\SIMPLE}}=\exp{\left(\frac{1-d^2}{4}-\frac{d^3}{12n}+O\left(\frac{
d^2}{n}\right)\right)}.
\end{equation}
\end{thm}
This estimate on the probability of {\SIMPLE} combined with
(\ref{e:ConfSimple}) will enable us to bound the probability of
events in our regular graph model $\Gnd$ based on bounds on the
probabilities of the corresponding events in the ``easier"
Configuration Model $\Pnd$.

We now introduce a few basic facts on $\Pnd$ which will be useful
for later computations on this model. Let $e_a$ and $e_b$ be two
distinct elements of our set of $dn$ elements. We define the
indicator variable $I_{\{e_a,e_b\}}$ over $\Pnd$ for the event that
the pair $\{e_a,e_b\}$ is part of the random pairing. By the
symmetry on the model, $e_a$ is equally likely to be paired with any
other element, thus $\Prob{I_{\{e_a,e_b\}}=1}=\frac{1}{dn-1}$.

For any subset of indices $I\subseteq [n]$ we define the set
$T_I=\bigcup_{i\in I}c_{i}$. Fix two such disjoint subsets, $I$ and
$J$, where $|I|=t$ and $|J|=s$. Given a random pairing $P$ in
$\Pnd$, let $X_{T_I}(P)$ be the random variable counting the number
of pairs in $P$ that use only elements from $T_I$, and let
$X_{T_I,T_J}(P)$ be the random variable counting the number of pairs
in $P$ that use one element from $T_I$ and the other from $T_J$. By
linearity of expectation we have that
\begin{eqnarray}
\Exp{X_{T_I}}&=&\sum_{\{e_a,e_b\}\in{{T_I}\choose
2}}\Exp{I_{\{e_a,e_b\}}}={{dt}\choose 2}\frac{1}{dn-1}\label{e:ExpTI};\\
\Exp{X_{T_I,T_J}}&=&\sum_{(e_a,e_b)\in T_I\times
T_J}\Exp{I_{\{e_a,e_b\}}}=\frac{d^2st}{dn-1}.
\end{eqnarray}
The expected value of $X_{T_I,T_J}$ conditioned on the event
$X_{T_I}=i$ can be also be computed using the linearity of
expectation. There are $(dt-2i)ds$ potential pairs, and the
probability of each of these pairs to be in a random pairing is
$\frac{1}{dn-dt}$, thus
\begin{equation}\label{e:cexpTITJi}
\cExp{}{X_{T_I,T_J}}{X_{T_I}=i}=\frac{ds(dt-2i)}{dn-dt}.
\end{equation}

Let $P$ and $P'$ be two distinct pairings in $\Pnd$. We write
\begin{equation}\label{e:defswithcdif}
P\sim P'\;\Leftrightarrow\; \exists e_a,e_b\quad
P'=P\setminus\{\{e_a,e^P_a\},\{e_b,e^P_b\}\}\cup\{\{e_a,e_b\},\{e_a^P,e_b^P\}\},
\end{equation}
that is, $P\sim P'$ if $P$ and $P'$ differ only by a single simple
\emph{switch}. The following is a well known concentration result in the
Configuration Model which makes use of martingales, and the
Azuma-Hoeffding inequality (see e.g \cite{AlonSpencer2000},
\cite{McDiarmid98}).
\begin{thm}[\cite{Wormald99} Theorem 2.19]\label{t:RegAzuma}
If $X$ is a random variable defined on $\Pnd$ such that
$|X(P)-X(P')|\leq c$ whenever $P\sim P'$, then for all $\lambda>0$
\begin{eqnarray*}
\Prob{X\geq
\Exp{X}+\lambda}\leq\exp\left(\frac{-\lambda^2}{dnc^2}\right);\\
\Prob{X\leq
\Exp{X}-\lambda}\leq\exp\left(\frac{-\lambda^2}{dnc^2}\right).
\end{eqnarray*}
\end{thm}

A direct corollary of Theorem \ref{t:RegAzuma} and
(\ref{e:ConfSimple}) gives us the following concentration result for
$\Gnd$.
\begin{cor}[\cite{Wormald99}]\label{c:Azuma4RegGraphs} Let $Y$ be a random variable defined on
$\Gnd$ such that $Y(G(P))=X(P)$ for all $P\in{\SIMPLE}$ where $X$ is
a random variable defined on $\Pnd$ that satisfies the conditions of
Theorem \ref{t:RegAzuma}. Then for all $\lambda>0$
\begin{eqnarray*}
\Prob{{Y\geq \Exp{Y}+\lambda}}&\leq&
\frac{\exp{\left(\frac{-\lambda^2}{dnc^2}\right)}}{\Prob{{\SIMPLE}}};\\
\Prob{{Y\leq \Exp{Y}-\lambda}}&\leq&
\frac{\exp{\left(\frac{-\lambda^2}{dnc^2}\right)}}{\Prob{{\SIMPLE}}}.
\end{eqnarray*}
\end{cor}

\subsection{Working directly with $\Gnd$ - the Switching Technique}
A more recent approach, introduced by McKay in \cite{McKay81}, that
has come to be known as the \emph{Switching Technique}, enables us
to work directly on the random regular graph model, $\Gnd$, without
passing through the Configuration model, $\Pnd$, which becomes
futile for large values of  $d(n)$. This technique enables us to
overcome the basic difficulty of counting elements in $\Gnd$ by
giving an alternative ``relative" counting technique. The basic
operation is the following:
\begin{defn}
Let $G$ be a $d$-regular graph, and
$S=(v_0,\ldots,v_{2r-1})\subseteq V(G)$ be some ordered set of $2r$
vertices of $G$ such that for any $1\leq i\leq r$
$\{v_{2i},v_{2i+1}\}\in E(G)$ and $\{v_{2i+1},v_{2i+2}\}\notin E(G)$
(where the addition in the indices is done modulo $2r$). A
\emph{r-switch} of $G$ by $S$ is the removal of all $r$ edges
$\{v_{2i},v_{2i+1}\}$ and the addition of the $r$ non-edges
$\{v_{2i+1},v_{2i+2}\}$ to $G$ as edges.
\end{defn}
The result of applying an $r$-switch operation on a $d$-regular
graph is a $d$-regular graph, and this $r$-switch operation is
obviously reversible. Now, let $Q$ be some integer-valued graph
parameter, and denote by $\mathcal{Q}_k$ the set of all graphs for
which $Q(G)=k$. We can now bound the ratio
$\frac{|\mathcal{Q}_k|}{|\mathcal{Q}_{k+1}|}$ by bounding the ratio
of the number of $r$-switch operations that take us from a graph
$G'$ where $Q(G')=k+1$ to a graph $G$ where $Q(G)=k$ and the number
of $r$-switch operations that take us from a graph $G$ where
$Q(G)=k$ to a graph $G'$ where $Q(G')=k+1$ for any integer $r$. For
a more detailed explanation, the interested reader is referred to
the proofs of Lemma \ref{l:switch_oneset} and Corollary \ref{c:AK3},
or to \cite{KrivelevichEtAl2001}, \cite{CooperEtAl2002a},
\cite{KimSudVu2002}, and \cite{KimSudVu2007} where the Switching
Technique is used extensively to prove various results on the random
regular graph model. Throughout this paper, we will only make use of
the $2$-switch operation, but there are some cases (e.g.
\cite{KrivelevichEtAl2001}, \cite{CooperEtAl2002a}) where the more
involved $3$-switch was used to overcome technical difficulties.

\section{Some structural properties of $\Gnd$}\label{s:structureGnd}
We proceed to prove a series of standard probabilistic claims about
the structure of a typical graph sampled from $\Gnd$ which will be
needed in the course of the proof of Theorem \ref{t:main2}.
\subsection{Edges spanned by a subset of vertices}\label{s:oneset} In
this section we will analyze the number of edges spanned by a set of
vertices. Our main motivation for this is to be able to prove some
technical lemmas on the distribution of edges spanned by subsets of
these cardinalities, in the spirit very similar to lemmas proved in
\cite{AloKri97} for the model $\Gnp$, that will be used in the
course of the proof of Theorem \ref{t:main2} in the succeeding
section.

The \emph{adjacency matrix} of a $d$-regular graph $G$ on $n$
vertices labeled by $[n]$, is the $n\times n$ binary matrix, $A(G)$,
where $A(G)_{ij}=1$ iff $(i,j)\in E(G)$. As $A(G)$ is real and
symmetric it has an orthogonal basis of real eigenvectors and all
its eigenvalues are real. We denote the eigenvalues of $A(G)$ in
descending order by $\lambda_1\geq\lambda_2 \ldots\geq \lambda_n$,
where $\lambda_1=d$ and its corresponding eigenvector is $j_n$ (the
$n\times 1$ all ones vector). Finally, let
$\lambda=\lambda(G)=\max\{|\lambda_2(G)|,|\lambda_n(G)|\}$, and call
such a graph $G$ an $(n,d,\lambda)$-graph. For an extensive survey
of fascinating properties of $(n,d,\lambda)$-graphs the reader is
referred to \cite{KriSud2006}. The celebrated expander mixing lemma
(see e.g. \cite{AlonSpencer2000} or \cite{Chung2004}) states roughly
that the smaller $\lambda$ is, the more random-like is the graph.
Here we present a simple variant of this lemma, which bounds the
number of edges spanned by any subset of vertices in an
$(n,d,\lambda)$-graph.
\begin{lem}[The Expander Mixing Lemma - Corollary 9.2.6 in \cite{AlonSpencer2000}]\label{l:expmixlem}
For every $(n,d,\lambda)$-graph $G=(V,E)$, every subset of vertices
$U\subseteq V$ satisfies
\begin{equation*}
\left|e(U)-{|U|\choose 2}\frac{d}{n}\right|\leq\lambda|U|.
\end{equation*}
\end{lem}

Broder et al. \cite{BroEtAl99}, extending a previous result of
Friedman, Kahn and Szemer{\'e}di \cite{FriKahSze89}, who used the
so-called trace method, estimate the ``typical'' second eigenvalue
of the $\Gnd$ for $d=o(\sqrt n)$.
\begin{thm}[Broder et al. \cite{BroEtAl99}]\label{t:BroEtAl99}
For $d=o(\sqrt n)$, w.h.p. $\lambda(\Gnd)=O(\sqrt d)$.
\end{thm}
With Lemma \ref{l:expmixlem} and Theorem \ref{t:BroEtAl99} at hand,
the following theorem is an immediate consequence.
\begin{thm}\label{t:oneset}
For every $d=o(\sqrt n)$ if $G=(V,E)$ is sampled from $\Gnd$ then
w.h.p. every subset of vertices $U\subseteq V$ satisfies
\begin{equation}\label{e:oneset}
\left|e(U)-{|U|\choose 2}\frac{d}{n}\right|=O(|U|\sqrt{d}).
\end{equation}
\end{thm}

The authors in \cite{Ben2005} give an alternative proof to Theorem
\ref{t:oneset} based purely on combinatorial techniques and not
relying on the spectral properties of the random graph $\Gnd$.

Although Theorem \ref{t:oneset} bounds the number of edges spanned
by a single set of vertices, the bound obtained does not meet our
needs to prove Theorem \ref{t:main2}, and we will need tighter
bounds. Therefore, we will start by analyzing the distribution of
the number of pairs spanned by a single set of indices in $\Pnd$.

Let $I$ be a subset of indices from $[n]$ of cardinality $|I|=t$.
Denote by $\mC_{T_I,k}$ be the set of all pairings in $\Pnd$
where there are exactly $k$ pairs that use only elements from $T_I$.
The possible values for $k$ are
\begin{equation}\label{e:kvals_oneset}
\max\left\{0,dt-\frac{dn}{2}\right\}\leq k\leq\frac{dt}{2}.
\end{equation}

To compute the cardinality of the set $\mC_{T_I,k}$ we count
all possible ways to choose $2k$ elements from $T_I$ and pair them
up, pair the rest of the elements from $T_I$ to elements outside
$T_I$, and pair the rest of the elements. It follows that
\begin{equation*}
|\mC_{T_Ik}|={dt\choose{2k}}\frac{(2k)!}{k!2^k}{{dn-dt}\choose{dt-2k}}
(dt-2k)!
\frac{(dn-2dt+2k)!}{\left(\frac{dn}{2}-dt+k\right)!2^{\frac{dn}{2}-dt+k}}.
\end{equation*}
Now, set
\begin{equation}\label{e:fdef}
f(k)=\frac{|\mC_{T_I,k+1}|}{|\mC_{T_I,k}|}=
\frac{1}{4}\cdot\frac{1}{k+1}\cdot\frac{(dt-2k)(dt-2k+1)}{\frac{dn}{2}-dt+k+1},
\end{equation}
and for convenience extend $f(k)$ to be a real valued function. Note
that if $k_1>k_2$, both values in the range of $k$ then
\begin{equation}\label{e:fratio}
\frac{f(k_1)}{f(k_2)}=
\frac{k_2+1}{k_1+1}\cdot\frac{dt-2k_1}{dt-2k_2}\cdot\frac{dt-2k_1+1}{dt-2k_2+1}\cdot
\frac{\frac{dn}{2}-dt+k_2+1}{\frac{dn}{2}-dt+k_1+1}\leq\frac{k_2+1}{k_1+1}.
\end{equation}
\begin{lem}\label{l:confmodel}
For any subset $I\subseteq [n]$ where $|I|=t$,
$|\mC_{T_I,k}|$ is monotonically increasing from
$k=\max\left\{0,dt-\frac{dn}{2}\right\}$ to
$\lfloor\Exp{X_{T_I}}\rfloor$, and monotonically decreasing from
$k=\lceil\Exp{X_{T_I}}\rceil$ to $k=\frac{dt}{2}$.
\end{lem}
\begin{proof}
Set $k_0=\Exp{X_{T_I}}$, and let $k'$ satisfy $f(k')=1$. Trivially,
$f(k)$ is monotonically decreasing when $k$ ranges from
$\max\left\{0,dt-\frac{dn}{2}\right\}$ to $\frac{dt}{2}$, therefore
it is enough to show $k_0-1\leq k'\leq k_0$. By solving for $k'$ we
get
\begin{equation*}
k'=\frac{(dt)^2-2dn+3t-4}{2dn+6}=\frac{dt(dt+3)}{2(dn+3)}-\frac{dn+2}{dn+3}.
\end{equation*}
Applying (\ref{e:ExpTI}) we have
\begin{equation*}
k_0-k'=\frac{dt(dt-1)}{2(dn-1)}-\frac{dt(dt+3)-2(dn+2)}{2(dn+3)}=
\frac{2dt(dt-dn)+2(dn+2)(dn-1)}{(dn-1)(dn+3)}.
\end{equation*}
The above function gets its minimal value for $t=\frac{n}{2}$ thus
$k_0-k'\geq\frac{\frac{3}{2}(dn)^2+2dn-4}{(dn-1)(dn+3)}\geq 0$. On
the other hand, again, by applying (\ref{e:ExpTI}),
\begin{equation*}
k'-k_0+1=\frac{dt(dt+3)}{2(dn+3)}-\frac{dn+2}{dn+3}-\frac{dt(dt-1)}{2(dn-1)}+1\geq
\frac{2dt(\frac{dn}{2}-dt)}{(dn-1)(dn+3)}\geq 0,
\end{equation*}
completing the proof of our claim.
\end{proof}

\begin{lem}\label{l:confmod_oneset}
Let $P\in \Pnd$ and $I\subseteq [n]$ be a fixed set of
indices. Then,\\
{\rm(a)} for any value of $\Delta$,
$\Pr\left[\left|X_{T_I}-\Exp{X_{T_I}}\right|\geq\Delta\right]\leq
\frac{d|I|}{2}e^{-\frac{\Delta^2}{4\Exp{X_{T_I}}+2\Delta+4}}$;\\
{\rm(b)} if $\Delta>\Exp{X_{T_I}}$ then
$\Pr\left[\left|X_{T_I}-\Exp{X_{T_I}}\right|\geq\Delta\right]\leq\frac{d|I|}{2}\cdot
e^{-\frac{\Delta}{2}\ln\frac{2\Exp{X_{T_I}}+\Delta+2}{2\Exp{X_{T_I}}+2}}$.
\end{lem}
\begin{proof}
Set $|I|=t$, $k_0=\Exp{X_{T_I}}={{dt}\choose 2}\frac{1}{dn-1}$,
$k_1=k_0+\Delta$, and $k_2=k_0-\Delta$. Following Lemma
\ref{l:confmodel} we know that $f(k_0-1)>1>f(k_0)$. We start by
proving claim {\rm(a)}.
\begin{equation}\label{e:ProbXTIk1}
\Prob{X_{T_I}=k_1}\leq
\frac{|\mC_{T_I,k_1}|}{|\mC_{T_I,k_0}|}=f(k_0)^{\Delta}\cdot\prod_{i=k_0}^{k_1-1}\frac{f(i)}{f(k_0)}\leq
\left(\frac{f\left(k_0+\frac{\Delta}{2}\right)}{f(k_0)}\right)^{\frac{\Delta}{2}}\leq
\left(\frac{k_0+1}{k_0+\frac{\Delta}{2}+1}\right)^{\frac{\Delta}{2}},
\end{equation}
where the last inequality follows from (\ref{e:fratio}). Similarly,
\begin{equation}\label{e:ProbXTIk2}
\Prob{X_{T_I}=k_2}\leq
\frac{|\mC_{T_I,k_2}|}{|\mC_{T_I,k_0}|}=
\left(\frac{1}{f(k_0)-1}\right)^{\Delta}\cdot\prod_{i=k_2}^{k_0-1}\frac{f(k_0-1)}{f(i)}\leq
\left(\frac{f\left(k_0-1\right)}{f(k_0-\frac{\Delta}{2}-1)}\right)^{\frac{\Delta}{2}}\leq
\left(\frac{k_0-\frac{\Delta}{2}}{k_0}\right)^{\frac{\Delta}{2}},
\end{equation}
where, again, the last inequality follows from (\ref{e:fratio}). To
prove claim {\rm(a)} we note that $\Prob{X_{T_I}=k_1}\leq
e^{-\frac{\Delta^2}{4k_0+2\Delta+4}}$. By Lemma \ref{l:confmodel} we
know that $|\mC_{T_I,k}|$ is monotonically increasing for
$k=\max\left\{0,dt-\frac{dn}{2}\right\}$ to $k_0$, and monotonically
decreasing for  $k=k_0$ to $\frac{dt}{2}$, therefore, by
(\ref{e:ProbXTIk1}) and (\ref{e:ProbXTIk2}), we have
\begin{eqnarray*}
&&\Pr\left[\left|X_{T_I}-k_0\right|\geq
\Delta\right]=\Prob{X_{T_I}\geq k_1}+\Prob{X_{T_I}\leq k_2}\leq\\
&&\left(\frac{dt}{2}-k_1\right)\cdot
e^{-\frac{\Delta^2}{4k_0+2\Delta+4}}+ k_2\cdot
e^{-\frac{\Delta^2}{4k_0}}\leq\frac{dt}{2}e^{-\frac{\Delta^2}{4k_0+2\Delta+4}}.
\end{eqnarray*}
The proof of claim {\rm(b)} is very similar. We note that
$\Prob{X_{T_I}=k_1}\leq
e^{-\frac{\Delta}{2}\ln\frac{2k_0+\Delta+2}{2k_0+2}}$ and that
$\Prob{X_{T_I}=k_2}\leq
e^{-\frac{\Delta}{2}\ln\frac{2k_0}{2k_0-\Delta}}$. By Lemma
\ref{l:confmodel} we know that $|\mC_{T_I,k}|$ is
monotonically increasing for
$k=\max\left\{0,dt-\frac{dn}{2}\right\}$ to $k_0$, and monotonically
decreasing for $k=k_0$ to $\frac{dt}{2}$, therefore, by
(\ref{e:ProbXTIk1}) and (\ref{e:ProbXTIk2}), we have
\begin{equation*}
\Prob{\left|X_{T_I}-k_0\right|\geq \Delta}=\Prob{X_{T_I}\geq
k_1}\leq\frac{dt}{2}\cdot
e^{-\frac{\Delta}{2}\ln\frac{2k_0+\Delta+2}{2k_0+2}}.
\end{equation*}
\end{proof}

We now apply the Switching Technique to analyze the distribution of
edges in a single subset of vertices. Fix a set $U\subseteq V(G)$ of
size $u$. Let $\mC_i$ denote the set of all $d$-regular
graphs where exactly $i$ edges have both ends in $U$.
\begin{lem}\label{l:switch_oneset}
For every $d,u=o(n)$,
\begin{equation*}
\frac{\left|\mC_i\right|}{\left|\mC_{i-1}\right|}\leq
\frac{1}{i}{u\choose 2}\frac{d}{n} \left(1+\frac{2u+4d}{n}\right).
\end{equation*}
\end{lem}
\begin{proof}
Let $F_i$ be a bipartite graph whose vertex set is composed of
$\mC_{i}\cup \mC_{i-1}$ and $\{G_1,G_2\}\in E(F_i)$
if and only if $G_1\in \mC_{i}$, $G_2\in \mC_{i-1}$,
and $G_2$ can be derived from $G_1$ by a 2-switch (and vice versa).
Let $G\in\mC_{i-1}$. To go to a graph in $\mC_{i}$
with a single $2$-switch operation, we need to choose the vertices
$a,b\in U$, where $\{a,b\}\notin E(G)$, and $x,y\in V(G)\setminus U$
such that $\{a,x\},\{b,y\}\in E(G)$ and $\{x,y\}\notin E(G)$, and
perform the switch as follows:
$\{a,x\},\{b,y\}\rightarrow\{a,b\},\{x,y\}$. This yields the upper
bound $d_{F_i}(G)\leq{u\choose 2}d^2$.

In order to give a lower bound on $d_{F_i}(G')$, where
$G'\in\mC_{i}$, we note that any choice of $a,b\in U$ such
that $\{a,b\}\in E(G)$ and an edge $\{x,y\}\in E(G')$ such that
$x,y\notin U$ and $x,y\notin N_{G'}(\{a,b\})$, gives us two ways to
perform the switch: $\{a,b\},\{x,y\}\rightarrow\{a,x\},\{b,y\}$ or
$\{a,b\},\{x,y\}\rightarrow\{a,y\},\{b,x\}$ resulting in a graph in
$C_{i-1}$. We note that $|N_{G'}(U)|\leq ud$ and
$|N_{G'}(\{a,b\})|\leq 2d$, thus $d_{F_i}(G')\geq
2i\left(\frac{dn}{2}-ud-d^2\right)$.

Combining the upper and lower bounds gives us
\begin{equation*}
\left|\mC_{i}\right|\cdot2i\left(\frac{dn}{2}-ud-2d^2\right)\leq
\sum_{G'\in\mC_{i}}d_{F_i}(G')=e(F_i)=
\sum_{G\in\mC_{i-1}}d_{F_i}(G)\leq
\left|\mC_{i-1}\right|\cdot {u\choose2}d^2.
\end{equation*}
Using the fact that $\frac{dn}{2}-ud-2d^2>ud+2d^2>0$ we conclude
that
$$\frac{|\mC_{i}|}{|\mC_{i-1}|}\leq
\frac{1}{i}\cdot{u\choose2}\frac{d}{n}\cdot\left(1+\frac{ud+2d^2}{\frac{dn}{2}-ud-2d^2}\right)\leq
\frac{1}{i}\cdot{u\choose2}\frac{d}{n}\cdot\left(1+\frac{2u+4d}{n}\right).$$
\end{proof}

The previous lemma give us the necessary ingredient to prove the following concentration result for $e(U)$.

\begin{cor}\label{c:switch_oneset}
For every $d,u=o(n)$, $\Delta\geq{u\choose 2}\frac{d}{n}$ and fixed
set of $u$ vertices $U$ in $\Gnd$,\\
{\rm(a)}$\Prob{e(U)\geq \Exp{e(U)}+\Delta}\leq
\frac{du}{2}\cdot\exp\left(-\frac{\Delta^2}{4\Exp{e(U)}+2\Delta}+o(\Delta)\right)$;\\
{\rm(b)}$\Prob{e(U)\geq\Exp{e(U)}+\Delta}\leq\frac{du}{2}\cdot
\exp\left(-\frac{\Delta}{2}\ln\left(1+\frac{\Delta}{2\Exp{e(U)}}\right)+o(\Delta)\right)$.
\end{cor}
\begin{proof}
Set $k_0=\Exp{e(U)}={u\choose{2}}\frac{d}{n}$ and $k=k_0+\Delta$. By
Lemma \ref{l:switch_oneset} it follows that
$\frac{|\mC_{i}|}{|\mC_{i-1}|}\leq\frac{1}{i}{u\choose2}\frac{d}{n}\left(1+o(1)\right)$,
and that $\left|\mC_{i}\right|$ is monotonically decreasing for
$i\geq 2{u\choose 2}\frac{d}{n}$ and hence,

\begin{eqnarray}\label{e:switch_oneset}\nonumber
&&\Prob{e(U)\geq k}\leq
\sum_{j=k}^{\frac{du}{2}}\frac{\left|\mC_{j}\right|}{\left|\mC_{k_0}\right|}\leq
\frac{du}{2}\cdot\prod_{i=k_0+1}^{k}\frac{|\mC_{i}|}{|\mC_{i-1}|}=\\
&&\frac{du}{2}\cdot(1+o(1))^{\Delta}\prod_{i=1}^{\Delta}\frac{k_0}{k_0+i}\leq
\frac{du}{2}\cdot
e^{o(\Delta)}\cdot\left(\frac{k_0}{k_0+\frac{\Delta}{2}}\right)^{\frac{\Delta}{2}}.
\end{eqnarray}

On the one hand, the right term of (\ref{e:switch_oneset}) equals
$\frac{du}{2}\cdot\exp\left(-\frac{\Delta}{2}\ln\frac{2k_0+\Delta}{2k_0}+o(\Delta)\right)$,
proving claim {\rm(b)}. On the other hand, $\frac{du}{2}\cdot
e^{o(\Delta)}\cdot\left(\frac{k_0}{k_0+\frac{\Delta}{2}}\right)^{\frac{\Delta}{2}}\leq
\frac{du}{2}\cdot\exp\left(-\frac{\Delta^2}{4k_0+2\Delta}+o(\Delta)\right)$,
proving claim {\rm(a)}.
\end{proof}
The concentration of $e(U)$ provided by Corollary
\ref{c:switch_oneset} enables us to obtain some upper bounds on the
number of edges spanned by subsets of vertices of different
cardinalities.
\begin{cor}\label{c:oneset1}
For every constant $C>0$ and $d=o(n^{1/5})$, w.h.p. every subset of
$u\leq C\sqrt{nd^3}$ vertices of $\Gnd$ spans less than $5u$ edges.
\end{cor}
\begin{proof}
Let $C$ be some positive constant. For every vertex subset $U$ of
cardinality $u\leq C\sqrt{nd^3}=o(n^{4/5})$, let
$\Delta=\Delta(U)=4u$, and $k_0=\Exp{e(U)}={u\choose
2}\frac{d}{n}\ll \Delta$. Note that we can assume that $u>5$, since
for $u\leq 5$ the claim is trivial. Summing over all possible values
of $u$, and applying Corollary \ref{c:switch_oneset} {\rm(b)} and
the union bound, we have that the probability that there exists a
set $U$ of at most $C\sqrt{nd^3}$ vertices, that spans more than
$4u+k_0<5u$ edges is bounded by
\begin{eqnarray*}
&&\sum_{u=6}^{C\sqrt{nd^3}}{n\choose{u}}\frac{du}{2}\cdot\exp\left(-2u\ln\left(1+\frac{4n}{(u-1)d}\right)+o(u)\right)\leq\\
&&\sum_{u=6}^{C\sqrt{nd^3}}\exp\left(-2u\ln\frac{4n}{ud}+u\ln\frac{n}{u}+u+\ln(du)+o(u)\right)\leq\\
&&n\cdot\exp\left(-5\ln n\right)=o(1).
\end{eqnarray*}
\end{proof}

\begin{cor}\label{c:oneset3}
For every $d=o(n)$ w.h.p. every subset of $u\geq\frac{n\ln n}{d}$
vertices of $\Gnd$ spans less than $\frac{5u^2d}{n}$ edges.
\end{cor}
\begin{proof}
For every vertex subset $U$ of cardinality $u\geq\frac{n\ln n}{d}$,
let $\Delta=\Delta(U)=\frac{4u^2d}{n}$, and
$k_0=\Exp{e(U)}={u\choose 2}\frac{d}{n}<\frac{u^2d}{2n}$. Summing
over all values of $u$, and applying Corollary \ref{c:switch_oneset}
{\rm(a)} and the union bound, we have that the probability that
there exists a set $U$ of at least $\frac{n\ln n}{d}$ vertices, that
spans more than $\frac{4u^2d}{n}+k_0<\frac{5u^2d}{n}$ vertices is
bounded by
\begin{eqnarray*}
&&\sum_{u=\frac{n\ln n}{d}}^{n}{n\choose{u}}\frac{du}{2}\cdot\exp\left(-\frac{8u^2d}{5n}+o\left(\frac{u^2d}{n}\right)\right)\leq\\
&&\sum_{u=\frac{n\ln
n}{d}}^{n}\exp\left(-\frac{3u^2d}{2n}+\ln(du)+u+u\ln\frac{n}{u}\right)\leq\\
&&\sum_{u=\frac{n\ln n}{d}}^{n}\exp\left(-\frac{3}{2}u\ln
n+\ln(du)+u+u\ln\frac{n}{u}\right)\leq\\
&&n\cdot\exp\left(-\frac{n\ln^2n}{3d}\right)=o(1).
\end{eqnarray*}
\end{proof}

\subsection{Other structural properties of $\Gnd$}
In the course of the proof of Theorem \ref{t:main2} we will need
some other asymptotic properties of $\Gnd$ which we proceed to
prove. We note that Lemma \ref{l:AK2} and Corollary \ref{c:AK3} can
also be obtained using a general result of Kim, Sudakov and Vu
\cite[Theorem 1.3]{KimSudVu2007} using density considerations, but
we prefer to give direct proofs using the Switching Technique,
similarly to what we have previously seen.

\begin{lem}\label{l:AK2}
For every $d=o(n^{1/5})$ w.h.p. every vertex participates in at most
$4$ triangles in $\Gnd$.
\end{lem}
\begin{proof}
Fix some vertex $v\in[n]$, and let $\mE_i$ denote the set of
all $d$-regular graphs on the set of vertices $[n]$ where
$e(N_G(v))=i$. We denote by $H_i$ the bipartite graph whose vertex
set is composed of $\mE_{i}\cup \mE_{i-1}$ and
$\{G_1,G_2\}\in E(H_i)$ if and only if $G_1\in \mE_{i}$,
$G_2\in \mE_{i-1}$, $G_2$ can be derived from $G_1$ by a
$2$-switch (and vice versa), and $N_{G_1}(v)=N_{G_2}(v)$ (that is,
we require that performing a $2$-switch does not affect the set of
$v$'s neighbors). This additional restriction enables us to apply
Lemma \ref{l:switch_oneset}\footnote{This statement requires some
caution. Recall that in Lemma \ref{l:switch_oneset}, if taking
$U$ as the set of neighbors of some vertex $v$, then counting the
number of 2-switches that take a graph from
$\mC_{i}$ to $\mC_{i-1}$ the neighbor set of $v$ may change. Nevertheless, the bound in the proof avoids
counting edges leaving $U$, so it applies to this case as
well.} to get that
$\frac{\left|\mE_i\right|}{\left|\mE_{i-1}\right|}\leq
\frac{1}{i}{d\choose
2}\frac{d}{n}\left(1+\frac{6d}{n}\right)<\frac{d^3}{in}=o(n^{-2/5})$.
It follows that for the given values of $d$, $|\mE_i|$ is
monotonically decreasing as $i$ goes from $0$ to ${d\choose 2}$. The
probability that there are more than four edges spanned by $N_G(v)$
can be bounded by
\begin{eqnarray*}
&&\Prob{e(N_G(v))\geq 5}\leq
\sum_{k=5}^{d\choose2}\frac{|\mE_{k}|}{|\mE_{0}|}\leq
\frac{|\mE_5|}{|\mE_0|}\cdot\sum_{j=0}^{{d\choose2}-5}\left(n^{-2/5}\right)^j=\\
&&\left(1+o(1)\right)\prod_{i=1}^{5}\frac{|\mE_{i}|}{|\mE_{i-1}|}\leq
\left(1+o(1)\right)\left(n^{-2/5}\right)^{5}=o(n^{-1}),
\end{eqnarray*}
and applying the union bound over all vertices from $[n]$ completes
the proof.
\end{proof}
Fix $U,W\subseteq V(G)$ to be two not necessarily disjoint subsets
of vertices of cardinality $u$ and $w$ respectively. Let
$\mD_{i}$ denote the subset of all $d$-regular graphs where
there are exactly $i$ edges with one endpoint in $U$ and the other
in $W$. The following proof follows closely the proof of Lemma \ref{l:switch_oneset}.
\begin{lem}\label{l:switch_2sets}
For $w\leq u\leq\frac{n}{10}$ and $d=o(n)$,
\begin{equation*}
\frac{|\mD_{i}|}{|\mD_{i-1}|}<\frac{2uwd}{in}.
\end{equation*}
Moreover, $|\mD_{i}|$ is monotonically decreasing for
$i>\frac{2uwd}{n}$.
\end{lem}
\begin{proof}
Let $H_i$ be a bipartite graph whose vertex set is composed of
$\mD_{i}\cup \mD_{i-1}$ and $\{G_1,G_2\}\in E(H_i)$
if and only if $G_1\in \mD_{i}$, $G_2\in \mD_{i-1}$,
and $G_2$ can be derived from $G_1$ by a $2$-switch (and vice
versa). Let $G\in\mD_{i-1}$. To give an upper bound on
$d_{H_i}(G)$ we need to find $a\in U$ and $b\in W$ such that
$\{a,b\}\notin E(G)$, and find $x\in V(G)\setminus W$ and $y\in
V(G)\setminus U$ such that $\{a,x\},\{b,y\}\in E(G)$ and
$\{x,y\}\notin E(G)$. Now we can perform the switch
$\{a,x\},\{b,y\}\rightarrow\{a,b\},\{x,y\}$ resulting in a graph in
$\mD_{i}$. This gives us an upper bound $d_{H_i}(G)\leq
uwd^2$.

On the other hand, let $G'\in\mD_{i}$. To derive a lower
bound on $d_{H_i}(G')$ we can first choose $a\in U$ and $b\in W$
such that $\{a,b\}\in E(G')$ in $i$ ways. Next we find an edge
$\{x,y\}\in E(G')$ such that $x,y\notin U\cup W$ and $x,y\notin
N_{G'}(\{a,b\})$. The number of edges with an endpoint in $U\cup W$
is at most $(u+w)d$ and $|N_{G'}(\{a,b\})|\leq 2d$. After choosing
such $x$ and $y$ we have two ways to perform the switch:
$\{a,b\},\{x,y\}\rightarrow\{a,x\},\{b,y\}$ or
$\{a,b\},\{x,y\}\rightarrow\{a,y\},\{b,x\}$ which gives us a graph
in $D_{i-1}$. It follows that
$d_{H_i}(G')\geq2i\left(\frac{dn}{2}-(u+w)d-2d^2\right)$.

We now compute the ratio, much like we have previously seen.
\begin{equation*}
\left|\mD_{i}\right|\cdot2i\left(\frac{dn}{2}-(u+w)d-2d^2\right)\leq
\sum_{G'\in\mD_{i}}d_{H_i}(G')=
\sum_{G\in\mD_{i-1}}d_{H_i}(G)\leq
\left|\mD_{i-1}\right|\cdot uwd^2.
\end{equation*}
Since $\frac{dn}{2}-(u+w)d-2d^2>(u+w)d+2d^2>0$, we have
$$\frac{|\mD_{i}|}{|\mD_{i-1}|}\leq
\frac{uwd}{in}\cdot\left(\frac{dn}{dn-2(u+w)-4d^2}\right)<\frac{2uwd}{in},$$
as claimed.
\end{proof}

\begin{cor}\label{c:AK3}
For every $d=o(n^{1/5})$ w.h.p. the number of paths of length three
between any two distinct vertices in $\Gnd$ is at most $4$.
\end{cor}
\begin{proof}
Fix $u,w\in[n]$ be two  distinct vertices. Let $\mF_i$
denote the set of all $d$-regular graphs, where there are exactly
$i$ paths of length three between $u$ and $w$. We denote by $H_i$
the bipartite graph whose vertex set is composed of
$\mF_{i}\cup \mF_{i-1}$ and $\{G_1,G_2\}\in E(H_i)$
if and only if $G_1\in \mF_{i}$, $G_2\in \mF_{i-1}$,
$G_2$ can be derived from $G_1$ by a $2$-switch (and vice versa),
$N_{G_1}(u)=N_{G_2}(u)$ and $N_{G_1}(w)=N_{G_2}(w)$ (that is,
performing the $2$-switch does not affect the neighbor sets of $u$
and $w$). We may therefore set $U=N_{G_1}(u)$ and $W=N_{G_1}(w)$ and
use Lemma \ref{l:switch_2sets}\footnote{Same as previous footnote.
The way we count the numbers of performing the 2-switch in Lemma
\ref{l:switch_2sets} will not affect the sets of neighbors of $u$
and $w$.} to derive that
$\frac{|\mF_{i}|}{|\mF_{i-1}|}<
\frac{2d^3}{in}=o(n^{-2/5})$. It follows, that for the given values
of $d$, $\mF_i$ is monotonically decreasing. Now, we bound
the probability that there are more than $5$ paths of length three
between $u$ and $w$ as follows:
\begin{eqnarray*}
&&\Prob{\hbox{There are more than }4\hbox{ paths of length
three between }u\hbox{ and }w}\leq\\
&&\sum_{k=5}^{n^2} \frac{|\mF_{k}|}{|\mF_{0}|}\leq\frac{|\mF_{5}|}{|\mF_{0}|}\cdot\sum_{j=0}^{n^2-5}\left(n^{-2/5}\right)^j\leq\\
&&\left(1+o(1)\right)\cdot\prod_{i=1}^{5}\frac{|\mF_{i}|}{|\mF_{i-1}|}\leq
\left(1+o(1)\right)\cdot o\left(n^{-2/5}\right)^{5}=o(n^{-2}).
\end{eqnarray*}
Using the union bound over all pairs of vertices from $[n]$
completes our proof.
\end{proof}

\section{Proof of Theorem \ref{t:main2}}\label{s:ChromConcentrarion}
In this section we prove Theorem \ref{t:main2}. The proof will
follow closely the proof of Alon and Krivelevich in \cite{AloKri97}
for the random graph model $\Gnp$.

Let us introduce the notion of graph choosability. A graph
$G=(\{v_1,\ldots,v_n\},E)$ is \emph{$\mathcal{S}$-choosable}, for a
family of color lists $\mathcal{S}=\{S_1,...,S_n\}$, if there exists
a proper coloring $f$ of $G$ that satisfies for every $1\leq i\leq
n$, $f(v_{i})\in S_i$. $G$ is \emph{$k$-choosable}, for some
positive integer $k$, if $G$ is $\mathcal{S}$-choosable for any
family $\mathcal{S}$ such that $|S_{i}|=k$ for every
$i\in\{1,...,n\}$. The \emph{choice number} of $G$, which is denoted
by $ch(G)$, is the minimum integer $k$ such that $G$ is
$k$-choosable.

A graph is \emph{$d$-degenerate} if every subgraph of it contains a
vertex of degree at most $d$. The following gives a trivial upper
bound on the choice number of a graph.
\begin{prop}\label{p:degchoose}
Every $d$-degenerate graph is $(d+1)$-choosable.
\end{prop}
For every $1\geq\varepsilon>0$ define $\tau=\tau(n,d,\varepsilon)$ to be
the least integer for which
\begin{equation}\label{e:taudef}
\Prob{\chi(\Gnd)\leq \tau}\geq\varepsilon
\end{equation}
and let $Y(G)$ be the random variable defined over $\Gnd$ that
denotes the minimal size of a set of vertices $S$ for which
$G\setminus S$ can be $\tau$-colored.
\begin{lem}\label{l:boundY(G)}
For every integer $n$,  $\varepsilon>0$ and $d=o(\sqrt n)$, there exists a constant,
$C=C(\varepsilon)$ such that
\begin{equation*}
\Prob{Y(\Gnd)\geq C\sqrt{nd^3}}\leq\varepsilon.
\end{equation*}
\end{lem}
\begin{proof}
Let $G$ be a random graph in $\Gnd$ and fix some $\varepsilon>0$. By
the minimality of $\tau$ it follows that
$\Prob{{\chi(G)<\tau}}<\varepsilon$. Define a proper coloring of the
multigraph generated by a pairing $P$, $G(P)$, as a proper coloring
of the this multigraph discarding its loops, and let $Y'$ be the
random variable defined over $\Pnd$ such that for every $P\in\Pnd$
$Y'(P)$ is the minimal size of a set of vertices $S'$ for which
$G(P)\setminus S'$ can be $\tau$-colored. Obviously, for
$P\in{\SIMPLE}$, $Y(G(P))=Y'(P)$. Let $P_0\in {\SIMPLE}$, and denote
by $P_0(m)$ the subset of pairs from $P_0$ which covers all of the
first $m$ elements. We define the following random variables over
$\Pnd$:
\begin{equation}\label{e:Y'martingale}
\forall 1\leq m\leq dn \qquad Y'_m(P_0)=\cExp{P\in
\Pnd}{Y'(P)}{P_0(m)\subseteq P},
\end{equation}
i.e., $Y'_m(P_0)$ is the expectation of the size of $S'$ conditioned
on all the pairings in $\Pnd$ that have the same first $m$ pairs as
$P_0$. $Y'_0(P_0),..., Y'_{dn-1}(P_0)$ is indeed a martingale, and
the random variable $Y'$ satisfies $|Y'(P)-Y'(P')|\leq 2$ for all
$P\sim P'$ (recalling \eqref{e:defswithcdif}). Setting
$\lambda=2\sqrt{d\ln\left(\varepsilon\cdot\Prob{{\SIMPLE}}\right)^{-1}}$,
and applying Corollary \ref{c:Azuma4RegGraphs} implies the following
concentration result on $Y$:
\begin{eqnarray}\nonumber\label{e:boundY(G)}
\Prob{Y(G)\geq
\Exp{Y(G)}+\lambda\sqrt{n}}\leq\frac{e^{-\lambda^2/4d}}{\Prob{{\SIMPLE}}}=\varepsilon;\\
\Prob{Y(G)\leq
\Exp{Y(G)}-\lambda\sqrt{n}}\leq\frac{e^{-\lambda^2/4d}}{\Prob{{\SIMPLE}}}=\varepsilon.
\end{eqnarray}
Notice that $\Prob{{Y(G)=0}}=\Prob{{\chi(G)\leq \tau}}>\varepsilon$,
therefore, by (\ref{e:boundY(G)}),
$\Exp{Y(G)}<\lambda\sqrt{n}$, and thus,
\begin{equation}\label{e:boundY(G)b}
\Prob{{Y(G)\geq 2\lambda\sqrt{n}}}\leq\Prob{{Y(G)\geq
\Exp{Y(G)}+\lambda\sqrt{n}}}\leq\varepsilon.
\end{equation}

Theorem \ref{t:BoundProbSimple} implies an upper bound on $\lambda$,
\begin{equation*}
\lambda=2\sqrt{-d(\ln{\varepsilon}+\ln{\Prob{{\SIMPLE}}})}\leq
2\sqrt{-d\left(\ln
\varepsilon+\left(\frac{1-d^2}{4}-\frac{d^3}{12n}-O\left(\frac{d^2}{n}\right)\right)\right)}
=O\left(\frac{d^{3/2}}{2}\right).
\end{equation*}
Returning to (\ref{e:boundY(G)b}), $\Prob{Y(G)\geq C\sqrt{nd^3}}\leq\varepsilon$ where $C$ is some constant depending on $\varepsilon$, as
claimed.
\end{proof}

Coming to prove Theorem \ref{t:main2}, we can use the previous
result of Achlioptas and Moore \cite{AchMoo2004}\footnote{In
\cite{AchMoo2004} the authors formally state this claim for the case
where $d$ is a constant (Theorem 1 in their paper), but in the first
paragraph of Section 2 of their paper they state it for
$d=O\left(n^{1/7-\delta}\right)$ for all $\delta>0$. It appears that
in their computations there may have been an oversight of the fact
that the probability of {\SIMPLE} needs to be taken into
consideration when moving from $\Pnd$ to $\Gnd$. Nevertheless, when
correcting this apparent oversight it can be shown that the proof
holds for $d=O\left(n^{1/9-\delta}\right)$ for all $\delta>0$ which
is the formulation we use here.}, whose proof of the same two-point
concentration result can be shown to go through for
$d=n^{1/9-\delta}$ for any $\delta>0$, and we can thus assume that
$d>n^{1/10}$. Let $\Gamma(n,d)$ denote the set of $d$-regular graphs
on $n$ vertices that satisfy the following properties:
\begin{enumerate}\label{Gamma}
\item {For every constant $C>0$, every subset of $u\leq C\sqrt{nd^3}$ vertices spans less than $5u$
edges.}\label{Gamma1}
\item {For every constant $C>0$, every subset of $u\leq Cn^{9/10}$ vertices spans less than $O(u\sqrt{d})$ edges.}\label{Gamma2}
\item {Every subset of $u\geq\frac{n\ln n}{d}$ vertices spans less than $O\left(\frac{u^2d}{n}\right)$ edges.}\label{Gamma3}
\item {For every vertex $v$, the number of edges spanned by $N(v)$ is at most $4$.}\label{Gamma4}
\item {The number of paths of length three between any two vertices is at most $4$.}\label{Gamma5}
\end{enumerate}

We have already proved that all properties of $\Gamma(n,d)$ occur
w.h.p. in $\Gnd$ for $d=o(n^{1/5})$ (Corollary \ref{c:oneset1},
Theorem \ref{t:oneset}, Corollary \ref{c:oneset3},  Lemma
\ref{l:AK2} and Corollary \ref{c:AK3} respectively), and by Theorem
\ref{t:FL92} we know that for $n^{1/10}<d\ll n^{1/5}$ w.h.p.
$\chi(\Gnd)\geq\frac{d}{2\ln d}>10$. The proof of Theorem
\ref{t:main2} now follows from the following deterministic
proposition, Proposition \ref{p:ChromConcentration}, by taking
$t=\tau(n,d,\frac{\varepsilon}{3})$, since by (\ref{e:taudef}) and
Lemma \ref{l:boundY(G)} we have:
\begin{eqnarray*}
&&\Prob{\chi(\Gnd)<t\hbox{ or
}\chi(\Gnd)>t+1}\leq\\
&&\Prob{\Gnd\notin\Gamma(n,d)}+ \frac{\varepsilon}{3}+\Prob{Y(\Gnd)\geq
c\sqrt{nd^3}}\leq \varepsilon.
\end{eqnarray*}
\begin{prop}\label{p:ChromConcentration}
There exists a positive integer $n_0$ such that for every $n\geq
n_0$, $n^{1/10}<d\ll n^{1/5}$, if $G\in\Gamma(n,d)$ such that
$\chi(G)\geq t\geq \frac{d}{2\ln d}$ and that there is a subset
$U_0\subseteq V$ of size $|U_0|=O(\sqrt{nd^3})$ such that
$G[V\setminus U_0]$ is $t$-colorable, then $G$ is $(t+1)$-colorable.
\end{prop}
\begin{proof}
First, we find a subset $U\subseteq V$ of size $O(\sqrt{nd^3})$
including $U_0$ such that every vertex $v\in V\setminus U$ has at
most $50$ neighbors in $U$. To find such a set $U$ we proceed as
follows. We start with $U=U_0$, and as long as there exists a vertex
$v\in V\setminus U$ with at least $50$ neighbors in $U$, we add $v$
to $U$ and iterate the process again. After $r$ iterations of this
process we have $|U|\leq c\sqrt{nd^3}+r$ and $e(U)\geq50r$. It
follows that the number of iterations is at most
$\frac{c\sqrt{nd^3}}{9}$, since otherwise, we would get a set of
$\frac{10c\sqrt{nd^3}}{9}$ vertices spanning at least
$\frac{50c\sqrt{nd^3}}{9}$ edges, a contradiction of Property
\ref{Gamma1} of the set $\Gamma(n,d)$. Let $U=\{u_1,\ldots,u_k\}$ be
the set at the end of the process, with $k=O(\sqrt{nd^3})$. Since by
Property \ref{Gamma1} of $\Gamma(n,d)$ every subset of
$i=O(\sqrt{nd^3})$ vertices spans less than $5i$ edges, then for
every $U'\subseteq U$ there is a vertex $v\in U'$ with
$d_{G[U']}(v)\leq\frac{2e(G[U'])}{|U'|}<10$. $G[U]$ is thus
$9$-degenerate and by Proposition \ref{p:degchoose}, $10$-choosable.

Let $f:V\setminus U\rightarrow \{1\ldots,t\}$ be a fixed proper
$t$-coloring of the subgraph $G[V\setminus U]$. Given this coloring
of $G[V\setminus U]$, we show that there exists a choice of $10$
color classes in each neighborhood $N_G(u_i)$ for every
$i\in\{1,\ldots,k\}$, such that the union of these $10k$ color
classes is an independent set. We can recolor the $10$ color classes
in each $N_G(u_i)$ by a new color $t+1$, yielding a proper
$(t+1)$-coloring, $g:V\setminus U\rightarrow \{1\ldots,t+1\}$, of
the vertices of $G[V\setminus U]$, and making $10$ colors from
$\{1,\ldots,t\}$ available for $u_i$. Since $G[U]$ is
$10$-choosable, there exists a proper coloring of $G[U]$ which
colors each $u_i$ from the set of $10$ colors available for it, that
extends $g$ to all of $G$, thus proving $G$ is $(t+1)$-colorable as
claimed.

We define an auxiliary graph $H=(W,F)$ whose vertex set is a
disjoint union of $k$ sets, $W_1,\ldots, W_k$, where for each vertex
$x\in N_G(U)$ and each neighborhood $N_G(u_i)$ in which it
participates there is a vertex in $w_{x,i}\in W_i$ corresponding to
$x$, and thus $|W_i|=|N_G(u_i)|=d$. We define
$\{w_{x,i},w_{y,j}\}\in F\Leftrightarrow \{x,y\}\in E(G[N_G(U)])$,
that is for every edge $\{x,y\}$ spanned in the edge set of
$N_G(U)$, we define the corresponding edges between all copies of
$x$ and $y$ in $H$. Since every $x\in N_G(U)$ has at most $50$
neighbors in $U$, there are at most $50$ copies of $x$ in $H$, and
thus every edge in $E(G[N_G(U)])$ yields at most $2500$ edges in
$H$. Furthermore, every independent set in $H$ corresponds to an
independent set in $G[N_G(U)]$, therefore, $f$ induces a proper
$t$-coloring, $f':W\rightarrow\{1,\ldots,t\}$, of the vertices of
$H$.

For every $s\leq k$ subsets $W_{i_1},\ldots,W_{i_s}$, the union has
$m=sd$ vertices. If $d\leq n^{4/25}$, then $m=O(n^{9/10})$ and thus,
by Property \ref{Gamma2} of $\Gamma(n,d)$, this union of sets spans
at most $O(m\sqrt d)$ edges in $H$. Now, for $n^{4/25}<d\ll
n^{1/5}$, if $s\leq n^{7/10}$ then $sd<n^{9/10}$ and, again by
Property \ref{Gamma2} of $\Gamma(n,d)$, this union spans at most
$O(m\sqrt d)$ edges, and if $s>n^{7/10}$ then $sd>\frac{n\ln n}{d}$,
and hence spans at most $O\left(\frac{m^2d}{n}\right)$ edges in $H$
by Property \ref{Gamma3} of $\Gamma(n,d)$. It follows, that for
every $s\leq k$ subsets $W_{i_1},\ldots,W_{i_s}$ there exists a set
$W_{i_l}$ connected by at most $O(m\sqrt{d}/s)=O(d^{3/2})$ edges to
the rest of the subsets if $d\leq n^{4/25}$ or $n^{4/25}<d\ll
n^{1/5}$ and $s\leq n^{7/10}$, and that there exists a set $W_{i_l}$
connected by at most $O\left(\frac{m^2d}{sn}\right)=
O\left(\frac{d^{9/2}}{n^{1/2}}\right)$ to the rest of the subsets if
$n^{4/25}<d\ll n^{1/5}$ and $s> n^{7/10}$. This implies that if
$d\leq  n^{4/25}$ the vertices $u_1,\ldots,u_k$ can be reordered in
such a way that for every $1<i\leq k$ there are $O(d^{3/2})$ edges
from $W_i$ to $\cup_{i'<i}W_{i'}$, and if $n^{4/25}<d\ll n^{1/5}$
the vertices $u_1,\ldots,u_k$ can be reordered in such a way that or
every $1<i\leq n^{7/10}$ there are $O(d^{3/2})$ edges from $W_i$ to
$\cup_{i'<i}W_{i'}$, and for every $n^{7/10}<i\leq k$ there are
$O\left(\frac{d^{9/2}}{n^{1/2}}\right)$ edges from $W_i$ to
$\cup_{i'<i}W_{i'}$. Assume that the vertices of $U$ are ordered in
such a way.

Now, according to the given order, we choose for each $u_i$, for $i$
ranging from $1$ to $k$, a set $J_i$ of $14$ colors. We say that a
color $c\in\{1,\ldots,t\}$ is \emph{available} for $u_i$ if there
does not exist an edge $\{w_{x,i},w_{y,i'}\}$ for some $i'<i$ such
that $f'(w_{x,i})=c$ and $f'(w_{y,i'})\in J_{i'}$, i.e., a color $c$
is available for $u_i$, if the corresponding color class in $W_i$,
has no connection with color classes having been chosen for previous
indices. The color lists $\{J_i\}_{i=1}^k$ are sequentially chosen
uniformly at random from the set of \emph{available} colors for
$u_i$.

Denote by $p_i$, for $1\leq i\leq k$, the probability that for some
$i'\leq i$ while choosing the set $J_{i'}$, there are less than
$\frac{t}{2}$ colors available for $i'$. Let us estimate $p_i$.
Obviously, $p_1=0$. First, assume $d\leq n^{4/25}$ and $1\leq i\leq
k$ or $n^{4/25}<d\ll n^{1/5}$ and $1<i\leq n^{7/10}$. In this case,
there are at most $O(d^{3/2})$ edges from $W_i$ to the previous sets
$W_{i'}$ for $i'<i$. Assuming $n^{4/25}<d\ll n^{1/5}$ and
$n^{7/10}<i\leq k$, it follows that there are at most
$O\left(\frac{d^{9/2}}{n^{1/2}}\right)$ edges from $W_i$ to the
previous sets $W_{i'}$ for $i'<i$. By Properties \ref{Gamma4} and
\ref{Gamma5} of $\Gamma(n,d)$, there are $\Theta(1)$ edges between
$W_{i'}$ and $W_i$, therefore each color chosen to be included in
$J_{i'}$ causes $\Theta(1)$ colors to become unavailable for $u_i$.
The probability of each color to be chosen into $J_{i'}$ is at most
$14$ divided by the number of available colors for $u_{i'}$ at the
moment of choosing $J_{i'}$. Hence, if $n^{1/10}< d\leq n^{4/25}$
and $1\leq i\leq k$ or $n^{4/25}<d\ll n^{1/5}$ and $1<i\leq
n^{7/10}$, then
\begin{eqnarray*}
&&p_i\leq
p_{i-1}+(1-p_{i-1}){{O(d^{3/2})}\choose{\frac{t/2}{\Theta(1)}}}
\left(\frac{14}{t/2}\right)^{\frac{t/2}{\Theta(1)}}\leq
p_{i-1}+\left(O(1)\frac{d^{3/2}}{t^2}\right)^{\Theta(t)}\leq\\
&&p_{i-1}+\exp\left(-C_1t\ln\frac{t^2}{d^{3/2}}\right)\leq
p_{i-1}+\exp\left(-\frac{C_1d}{2\ln
n}\ln\frac{d^{1/2}}{4\ln^2n}\right)\leq\\
&&p_{i-1}+e^{-C_2n^{1/10}},
\end{eqnarray*}
where $C_1$ and $C_2$ are positive constants. If $n^{4/25}<d\ll
n^{1/5}$ and $n^{7/10}<i\leq k$, then
\begin{eqnarray*}
&&p_i\leq
p_{i-1}+(1-p_{i-1}){{O\left(\frac{d^{9/2}}{n^{1/2}}\right)}\choose{\frac{t/2}{\Theta(1)}}}
\left(\frac{14}{t/2}\right)^{\frac{t/2}{\Theta(1)}}\leq
p_{i-1}+\left(O(1)\frac{d^{9/2}}{n^{1/2}t^2}\right)^{\Theta(t)}\leq\\
&&p_{i-1}+\exp\left(-C_3t\ln\frac{n^{1/2}t^2}{d^{9/2}}\right)\leq
p_{i-1}+\exp\left(-\frac{C_3d}{2\ln n}\ln\frac{n^{1/2}}{4d^{5/2}\ln^2n}\right)\leq\\
&&p_{i-1}+e^{-C_4n^{4/25}}\leq p_{i-1}+e^{-C_4n^{1/10}},
\end{eqnarray*}
where $C_3$ and $C_4$ are positive constants.

Since $p_k<O(ke^{-Cn^{1/10}})=o(1)$ for some constant $C$, there
exists a family of color lists $\{J_i:1\leq i\leq k,\;|J_i|=14\}$
for which there are no edges between the corresponding color classes
of distinct subsets $W_{i'},W_i$. Once such a family is indeed
found, for every $1\leq i\leq k$, we go over all edges inside $W_i$,
and for every edge choose one color class that is incident with it,
and delete its corresponding color from $J_i$. By Property
\ref{Gamma4} of $\Gamma(n,d)$, each $W_i$ spans at most four edges,
therefore, we deleted at most four colors from each $J_i$ and thus,
we get a family $\{I_i:1\leq i\leq k\;,\;|I_i|\geq 10\}$, for which
the union $\cup_{i=1}^k\{w\in W_i\;:\;f'(w)\in I_i\}$ is an
independent set in $H$, completing the proof.
\end{proof}
\section{Concluding remarks and open problems}
In this paper we proved that for $d=o(n^{1/5})$ the chromatic number
of a random $d$-regular graph on $n$ vertices is w.h.p. concentrated
on two consecutive integers. We propose here further questions which
may be of interest to pursue, but will most likely need new ideas to
resolve.
\begin{itemize}
\item Alon and Krivelevich in \cite{AloKri97} noted, using a
continuity argument, that for $\Gnp$ the two-value concentration is
best possible for a general $p\leq n^{1/2-\varepsilon}$ where
$\varepsilon >0$. On the other hand, they showed there exists a
series of values of $p$ in this range for which $\chi(\Gnp)$ is in
fact concentrated on a single value. This may as well be the case
for $\Gnd$, but as $d$ must be an integer, the arguments of Alon and
Krivelevich cannot be applied trivially.
\item Theorem \ref{t:main2} does not give any evidence as to what
are the actual values on which $\chi(\Gnd)$ is concentrated.
Achlioptas and Moore in \cite{AchMoo2004}, following ideas of
Achlioptas and Naor \cite{AchNao2005}, showed that for a constant
$d$, w.h.p. $\chi(\Gnd)$ is concentrated on three consecutive
integers $\{k,k+1,k+2\}$ where $k=k(d)=\min\{t\in\mathbb{N}:2t\ln
t>d\}$. In addition if $d>(2k-1)\log k$ then $\chi(\Gnd)$ is
concentrated on two consecutive integers $\{k+1,k+2\}$. Recently,
Kemkes, P\'{e}rez-Gim\'{e}nez and Wormald \cite{KemPerWorPre} showed
that w.h.p. $\chi(\Gnd)<k+2$ (extending \cite{ShiWor2007a,
ShiWor2007b}; see also related result on $\chi(\mathcal{G}_{n,5})$
in \cite{DiazEtAlPre}) thus determining exactly or up to two
consecutive integers $\chi(\Gnd)$. Albeit the above, locating the
concentration interval for non-constant values of $d$ remains open.
\item Lastly, the range of $d$ for which Theorem \ref{t:main2}
holds, seems to be far from optimal. The main obstacle to increase
the range of $d$ so as to match the corresponding result of Alon and
Krivelevich (i.e. taking $d$ to be as high as $n^{1/2-\varepsilon}$
for any $\varepsilon>0$) is the fact the in $\Pnd$ we could not find
a ``vertex-exposure'' martingale which satisfied a Lipschitz
condition. By using an analogue of the ``edge-exposure'' martingale,
our concentration result was much more restrictive on $d$. This
seems more of a technicality of our proof approach, and we believe
that Theorem \ref{t:main2} holds for a larger range of $d$.
\end{itemize}
\subsection*{Acknowledgement}
We would like to thank Alan Frieze for bringing Theorem
\ref{t:BroEtAl99} to our attention, Nick Wormald and Xavier
P\'{e}rez-Gim\'{e}nez for providing us several of the references
below and the anonymous referees for helpful comments and
corrections.
\bibliographystyle{abbrv}
\bibliography{regrand}

\begin{thebibliography}{10}

\bibitem{AchMoo2004}
D.~Achlioptas and C.~Moore.
\newblock The chromatic number of random regular graphs.
\newblock In {\em RANDOM: International Workshop on Randomization and
  Approximation Techniques in Computer Science}, pages 219--228. LNCS, 2004.

\bibitem{AchNao2005}
D.~Achlioptas and A.~Naor.
\newblock The two possible values of the chromatic number of a random graph.
\newblock {\em Annals of Mathematics. Second Series}, 162(3):1335--1351, 2005.

\bibitem{AloKri97}
N.~Alon and M.~Krivelevich.
\newblock The concentration of the chromatic number of random graphs.
\newblock {\em Combinatorica}, 17(3):303--313, 1997.

\bibitem{AlonSpencer2000}
N.~Alon and J.~H. Spencer.
\newblock {\em The Probabilistic Method, 2nd edition}.
\newblock Wiley-Interscience Series in Discrete Mathematics and Optimization.
  John Wiley \& Sons, 2000.

\bibitem{Ben2005}
S.~Ben-Shimon.
\newblock On the distribution of edges in random regular graphs.
\newblock Master's thesis, Tel Aviv University, 2005.
\newblock http://cs.tau.ac.il/research/sonny.benshimon/papers/MScThesis.pdf.

\bibitem{BenCan78}
E.~A. Bender and E.~R. Canfield.
\newblock The asymptotic number of non-negative integer matrices with given row
  and column sums.
\newblock {\em Journal of Combinatorial Theory, Series A}, 24(3):296--307,
  1978.

\bibitem{Bol80}
B.~Bollob\'{a}s.
\newblock A probabilistic proof of an asymptotic formula for the number of
  labelled regular graphs.
\newblock {\em European Journal of Combinatorics}, 1:311--316, 1980.

\bibitem{Bol88}
B.~Bollob\'{a}s.
\newblock The chromatic number of random graphs.
\newblock {\em Combinatorica}, 8(1):49--55, 1988.

\bibitem{Bol2001}
B.~Bollob\'{a}s.
\newblock {\em Random Graphs}.
\newblock Cambridge University Press, 2001.

\bibitem{BroEtAl99}
A.~Broder, A.~Frieze, S.~Suen, and E.~Upfal.
\newblock Optimal construction of edge-disjoint paths in random graphs.
\newblock {\em {SIAM} Journal on Computing}, 28(2):541--573, 1999.

\bibitem{Chung2004}
F.~Chung.
\newblock Discrete isoperimetric inequalities.
\newblock In A.~Grigor'yan and S.~T. Yau, editors, {\em Surveys in Differential
  Geometry: Eigenvalues of Laplacians and other geometric operators},
  volume~IX. International Press, 2004.

\bibitem{CojPanSte2008}
A.~Coja-Oghlan, K.~Panagiotou, and A.~Steger.
\newblock On the chromatic number of random graphs.
\newblock {\em Journal of Combinatorial Theory, Series B}, 98(5):980--993,
  2008.

\bibitem{CooperEtAl2002a}
C.~Cooper, A.~Frieze, B.~Reed, and O.~Riordan.
\newblock Random regular graphs of non-constant degree: Independence and
  chromatic number.
\newblock {\em Combinatorics, Probability and Computing}, 11(4):323--342, 2002.

\bibitem{DiazEtAlPre}
J.~D\'{\i}az, A.~C. Kaporis, G.~D. Kemkes, L.~M. Kirousis, X.~P\'{e}rez, and
  N.~C. Wormald.
\newblock On the chromatic number of a random $5$-regular graph.
\newblock {\em Journal of Graph Theory}, To appear.

\bibitem{FriKahSze89}
J.~Friedman, J.~Kahn, and E.~Szemer{\'e}di.
\newblock On the second eigenvalue in random regular graphs.
\newblock In {\em Proceedings of the Twenty First Annual {ACM} Symposium on
  Theory of Computing}, pages 587--598, 1989.

\bibitem{FriezeLuczak92}
A.~M. Frieze and T.~{\L}uczak.
\newblock On the independence and chromatic numbers of random regular graphs.
\newblock {\em Journal of Combinatorial Theory, Series B}, 54(1):123--132,
  1992.

\bibitem{JanLucRuc2000}
S.~Janson, T.~{\L}uczak, and A.~Ruci{\'n}ski.
\newblock {\em Random Graphs}.
\newblock Wiley-Interscience Series in Discrete Mathematics and Optimization.
  John Wiley \& Sons, 2000.

\bibitem{KemPerWorPre}
G.~D. Kemkes, X.~P\'{e}rez-Gimen\'{e}z, and N.~C. Wormald.
\newblock On the chromatic number of $d$-regulat graphs.
\newblock http://arxiv.org/abs/0812.2973.

\bibitem{KimSudVu2002}
J.~H. Kim, B.~Sudakov, and V.~H. Vu.
\newblock On the asymmetry of random regular graphs and random graphs.
\newblock {\em Random Structures \& Algorithms}, 21(3-4):216--224, 2002.

\bibitem{KimSudVu2007}
J.~H. Kim, B.~Sudakov, and V.~H. Vu.
\newblock Small subgraphs of random regular graphs.
\newblock {\em Discrete Mathematics}, 307(15):1961--1967, 2007.

\bibitem{KimVu2004}
J.~H. Kim and V.~H. Vu.
\newblock Sandwiching random graphs: Universality between random graph models.
\newblock {\em Advances in Mathematics}, 188(2):444--469, 2004.

\bibitem{KriSud2006}
M.~Krivelevich and B.~Sudakov.
\newblock Pseudo-random graphs.
\newblock In E.~Gy\H{o}ri, G.~O.~H. Katona, and L.~Lov\'{a}sz, editors, {\em
  More sets, graphs and numbers: A Salute to Vera S\`{o}s and Andr\'{a}s
  Hajnal}, volume~15 of {\em Bolyai Society Mathematical Studies}, pages
  199--262. Springer, 2006.

\bibitem{KrivelevichEtAl2001}
M.~Krivelevich, B.~Sudakov, V.~H. Vu, and N.~C. Wormald.
\newblock Random regular graphs of high degree.
\newblock {\em Random Structures \& Algorithms}, 18(4):346--363, 2001.

\bibitem{Luc91b}
T.~{\L}uczak.
\newblock The chromatic number of random graphs.
\newblock {\em Combinatorica}, 11(1):45--54, 1991.

\bibitem{Luc91}
T.~{\L}uczak.
\newblock A note on the sharp concentration of the chromatic number of random
  graphs.
\newblock {\em Combinatorica}, 11(3):295--297, 1991.

\bibitem{McDiarmid98}
C.~McDiarmid.
\newblock Concentration.
\newblock In M.~Habib, C.~McDiarmid, J.~Ramirez-Alfonsin, and B.~Reed, editors,
  {\em Probabilistic Methods for Algorithmic Discrete Mathematics}, pages
  195--248. Springer, 1998.

\bibitem{McKay81}
B.~D. McKay.
\newblock Subgraphs of random graphs with specified degrees.
\newblock {\em Congressus Numerantium}, 33:213--223, 1981.

\bibitem{McKayWormald91}
B.~D. McKay and N.~C. Wormald.
\newblock Asymptotic enumeration by degree sequence of graphs with degrees
  $o(n^{1/2})$.
\newblock {\em Combinatorica}, 11(4):369--382, 1991.

\bibitem{ShaSpe87}
E.~Shamir and J.~Spencer.
\newblock Sharp concentration of the chromatic number on random graphs
  $\mathcal{G}(n,p)$.
\newblock {\em Combinatorica}, 7(1):131--138, 1987.

\bibitem{ShiWor2007a}
L.~Shi and N.~C. Wormald.
\newblock Colouring random 4-regular graphs.
\newblock {\em Combinatorics, Probability and Computing}, 16(2):309--344, 2007.

\bibitem{ShiWor2007b}
L.~Shi and N.~C. Wormald.
\newblock Colouring random regular graphs.
\newblock {\em Combinatorics, Probability and Computing}, 16(3):459--494, 2007.

\bibitem{Wormald99}
N.~C. Wormald.
\newblock Models of random regular graphs.
\newblock In J.~Lamb and D.~Preece, editors, {\em Surveys in Combinatorics},
  volume 276 of {\em London Mathematical Society Lecture Note Series}, pages
  239--298. Cambridge University Press, 1999.

\end{thebibliography}
\end{document}